\documentclass[12pt]{amsart}
\usepackage{graphicx}
\usepackage{amsfonts}
\usepackage{amssymb}
\usepackage{latexsym}
\usepackage{amscd}

\catcode`\@=12

\def\Empty{}
\newcommand\oplabel[1]{
  \def\OpArg{#1} \ifx \OpArg\Empty {} \else
  	\label{#1}
  \fi}
		
\newcommand{\comm}[1]{}

\input{psfig}
\newtheorem{thm}{Theorem}[section]
\newtheorem{cor}[thm]{Corollary}
\newtheorem{lem}[thm]{Lemma}
\newtheorem{prop}[thm]{Proposition}

\newtheorem{fsc}{Fatou-Sullivan Classification of Fatou Components}

\newenvironment{pf}{\proof[\proofname]}{\endproof}
\newenvironment{pf*}[1]{\proof[#1]}{\endproof}
\usepackage{euscript}

\usepackage[OT2,OT1]{fontenc}

\newcommand{\cal}[1]{{\mathcal #1}}

\newcommand{\beq}{\begin{equation}}
\newcommand{\eeq}{\end{equation}}
\newcommand{\eref}[1]{(\ref{#1})}

\newcommand{\ve}{\varepsilon}

\newcommand{\al}{\alpha}

\theoremstyle{definition}
\newtheorem{defn}{Definition}[section]

\theoremstyle{remark}

\renewcommand{\deg}{\operatorname{deg}}
\newcommand{\riem}{\hat{\CC}}

\newcommand{\dist}{\operatorname{dist}}

\newcommand{\cl}{\operatorname{cl}}

\newcommand{\eps}{\epsilon}

\numberwithin{equation}{section}
\newcommand{\thmref}[1]{Theorem~\ref{#1}}
\newcommand{\propref}[1]{Proposition~\ref{#1}}

\newcommand{\lemref}[1]{Lemma~\ref{#1}}

\newcommand{\cC}{{\cal C}}

\newcommand{\cR}{{\cal R}}

\newcommand{\CC}{{\mathbb C}}
\newcommand{\RR}{{\mathbb R}}

\newcommand{\ZZ}{{\mathbb Z}}
\newcommand{\NN}{{\mathbb N}}
\newcommand{\DD}{{\mathbb D}}

\newcommand{\QQ}{{\mathbb Q}}

\begin{document}
\addtolength{\evensidemargin}{-0.7in}
\addtolength{\oddsidemargin}{-0.7in}

\title[Computability of filled Julia sets]{Filled Julia sets with empty interior are computable}
\author{I. Binder, M. Braverman, M. Yampolsky}
\thanks{The first and third authors are partially supported by NSERC 
Discovery grants; the second author is partially supported by NSERC 
Postgraduate Scholarship}
\date{\today}

\begin{abstract}
We show that if a polynomial filled Julia set has empty interior, then it is computable.

\end{abstract}
\maketitle

\section{Introduction}

\subsection*{Julia sets of rational mappings}
We recall the main definitions of complex dynamics relevant to our results only briefly;
a good general reference is the book of Milnor \cite{Mil}.
For a rational mapping $R$ of degree $\deg R=d\geq 2$ considered as a dynamical system on the Riemann
sphere
$$R:\riem\to\riem$$
the Julia set is defined as the complement of the set where the dynamics is Lyapunov-stable:

\begin{defn}
Denote $F_R$ the set of points $z\in\riem$ having an open neighborhood $U(z)$ on which the
family of iterates $R^n|_{U(z)}$ is equicontinuous. The set $F_R$ is called the Fatou set of $R$
and its complement $J_R=\riem\setminus F_R$ is the Julia set.
\end{defn}

\noindent
In the case when the rational mapping is a polynomial $$P(z)=a_0+a_1z+\cdots+a_dz^d:\CC\to\CC$$ an equivalent
way of defining the Julia set is as follows. Obviously, there exists a neighborhood of $\infty$ on $\riem$
on which the iterates of $P$ uniformly converge to $\infty$. Denoting $A(\infty)$ the maximal such domain of attraction
of $\infty$ we have $A(\infty)\subset F(R)$. We then have 
$$J_P=\partial A(\infty).$$
The bounded set $\riem \setminus \cl A(\infty)$ is called {\it the filled Julia set}, and denoted $K_P$;
it consists of points whose orbits under $P$ remain bounded:
$$K_P=\{z\in\riem|\;\sup_n|P^n(z)|<\infty\}.$$

\noindent
For future reference, let us list in a proposition below the main properties of Julia sets:

\begin{prop}
\label{properties-Julia}
Let $R:\riem\to\riem$ be a rational function. Then the following properties hold:
\begin{itemize}
\item $J_R$ is a non-empty compact subset of $\riem$ which is completely
invariant: $R^{-1}(J(R))=J(R)$;
\item $J_R=J_{R^n}$ for all $n\in\NN$;
\item $J_R$ is perfect;
\item if $J_R$ has non-empty interior, then it is the whole of $\riem$;
\item let $U\subset\riem$ be any open set with $U\cap J_R\neq \emptyset$. Then there exists $n\in\NN$ such that
$R^n(U)\supset J_R$;
\item periodic orbits of $R$ are dense in $J_R$.
\end{itemize}
\end{prop}

\noindent
Let us further comment on the last property. For a periodic point $z_0=R^p(z_0)$
of period $p$ its {\it multiplier} is the quantity $\lambda=\lambda(z_0)=DR^p(z_0)$.
We may speak of the multiplier of a periodic cycle, as it is the same for all points
in the cycle by the Chain Rule. In the case when $|\lambda|\neq 1$, the dynamics
in a sufficiently small neighborhood of the cycle is governed by the Mean
Value Theorem: when $|\lambda|<1$, the cycle is {\it attracting} ({\it super-attracting}
if $\lambda=0$), if $|\lambda|>1$ it is {\it repelling}.
Both in the attracting and repelling cases, the dynamics can be locally linearized:
\begin{equation}
\label{linearization-equation}
\psi(R^p(z))=\lambda\cdot\psi(z)
\end{equation}
where $\psi$ is a conformal mapping of a small neighborhood of $z_0$ to a disk around $0$.
By a classical result of Fatou, a rational mapping has at most finitely many non-repelling
periodic orbits.
The sharp bound on the number of such orbits is $2d-2$; it is established by Shishikura.
 Therefore, we may refine the last statement of \propref{properties-Julia}:

\begin{itemize}
\item {\it repelling periodic orbits are dense in $J_R$}.
\end{itemize}

\noindent
In the case when $|\lambda|=1$, so that $\lambda=e^{2\pi i\theta}$, $\theta\in\RR$, 
 the simplest to study is the {\it parabolic case} when $\theta=n/m\in\QQ$, so $\lambda$ 
is a root of unity. In this case $R^p$ is not locally linearizable; it is not hard to see that $z_0\in J_R$.
 In the complementary situation, two non-vacuous possibilities  are considered:
{\it Cremer case}, when $R^p$ is not linearizable, and {\it Siegel case}, when it is.
In the latter case, the linearizing map $\psi$ from (\ref{linearization-equation}) conjugates
the dynamics of $R^p$ on a neighborhood $U(z_0)$ to the irrational rotation by angle $\theta$
(the {\it rotation angle})
on a disk around the origin. The maximal such {\it rotation domain} around $z_0$ is called a {\it Siegel disk}.

A different kind of a rotation domain may occur for a non-polynomial rational mapping $R$.
A {\it Herman ring} $A$ is a conformal image $$\nu:\{z\in\CC|\;0<r<|z|<1\}\to A,$$
such that $$R^p\circ\nu(z)=\nu(e^{2\pi i\theta}z),$$
for some $p\in\NN$ and $\theta\in\RR\setminus\QQ$. 

To conclude our discussion of the basic facts of rational dynamics, let us formulate the Fatou-Sullivan
Classification of the types of connected components of the Fatou set of a rational mapping. 
The term {\it basin} in what follows will describe the set of points whose orbits converge to a 
given attracting or parabolic periodic orbit under the iteration of $R$. 

\begin{fsc}
For every connected component $W\subset F_R$ there exists $m\in\NN$ such that the image $H=R^m(W)$ is
periodic under the dynamics of $R$. Moreover, each periodic Fatou component $H$ is of one of the following
types: 
\begin{itemize}
\item a component of the basin  of  an attracting or a super-attracting periodic orbit;
\item a component of the  basin  of a parabolic periodic orbit;
\item a Siegel disk;
\item a Herman ring.
\end{itemize}
\end{fsc}

\noindent
\subsection*{Computability of subsets of $\RR^n$}
The computability notions that we apply for subsets of $\RR^k$ belong to 
the framework of {\it computable analysis}. 
Their roots can be traced to the pioneering work of Banach and Mazur of 1937 (see \cite{Maz}).
The reader may find an extended exposition of 
the model of computation we are using in \cite{BC2}. 
See \cite{Wei} for a detailed discussion of the concepts of modern computable analysis.
See also \cite{BY} for a discussion of 
computability as applied to problems in Complex Dynamics.

Given a compact set $S\subset \RR^k$, our goal is to be able to approximate 
the set $S$ with an arbitrarily high precision. Here we ought to specify 
what do ``approximate" and ``precision" mean in this setting. 

Denote $B(\bar x,r)\subset \RR^k$ the closed Euclidean ball with radius $r$ centered at 
$\bar x$.  
Denote $\cC$ the set of finite unions $\overline{\cup B(\bar x_i,r_i)}\subset\RR^k$ of closed balls whose
radii and coordinates of the centers are all dyadic rationals:
$$
\cC = \left\{ \bigcup_{i=1}^k B(\bar x_i, r_i)~:~\bar x_i \in \QQ^k,~r_i\in \QQ\right\}. 
$$
For a set $K \subset \RR^k$ denote its $\ve$-neigborhood by $B(K,\ve)=\cup_{\bar x\in K} B(\bar x,\ve)$.
Recall that the {\em Hausdorff metric} on compact subsets of $\RR^k$ is defined by
$$
d_H(X,Y) = \inf \{ \ve~:~ X \subset B(Y,\ve)\mbox{ and } Y\subset B(X,\ve)\}.
$$
That is, $d_H(X,Y)$ is the smallest quantity by which we have to ``expand" $X$ to cover 
$Y$ and vice-versa. 

In computable analysis
a compact set $S \subset \RR^k$ is called {\it computable}, if there exists a Turing Machine (TM) 
$M(m)$, such that on input $m$, $M(m)$ outputs an encoding of  $C_m \in \cC$ satisfying $d_H (S, C_m) < 2^{-m}$. 
That is, $S$ is effectively approximable in the Hausdorff metric. This definition is quite 
robust. For example, it is equivalent to the computability of the distance function
$d_S(x) = \inf \{|x-y|~:~y\in S\}$. 
In more practical terms, a set is computable if there exists a computer program which
generates its image on a computer screen with one-pixel precision for an arbitrarily
fine pixel size.

In the case of computing Julia sets, we are not actually computing a single set, but 
a set-valued function that takes the coefficients of the rational function $R(z)$ as an input, and 
outputs approximations of $J_R$. The machine computing $J_R$ is given access to $R(z)$ through 
an {\em oracle} that can  provide an approximation of any coefficient of $R(z)$ with an arbitrarily 
high (but finite) requested precision. Such an oracle machine is denoted by $M^\phi$. If 
$P$ is the set of rational functions on which $M^\phi$ works correctly, then $M^\phi$ is
said to be {\em uniform} on $P$. In case $P$ is a singleton, the machine is said to be 
{\em non-uniform}. In other words, a uniform machine is designed to work on some family 
of inputs, whereas a non-uniform machine is only supposed to work on a specific input. 

We remark that Julia sets have appeared in a somewhat similar context in the book of 
Blum, Cucker, Shub, and Smale \cite{BCSS}. They considered the question whether
a Julia set is decidable in the the Blum-Shub-Smale (BSS) model of real computation.
The BSS model is very different from the computable analysis model we use, and
can be very roughly described as based on computation with infinite-precision real arithmetic.
Some discussion of the differences between the models may be found in \cite{BC2} and the
references therein.
It turns out (see Chapter 2.4 of \cite{BCSS}) that in the BSS model all, 
but the most trivial Julia sets 
are not decidable. More generally, sets with a fractional Hausdorff dimension, 
including the ones with very simple description, such as the Cantor set, are 
BSS-undecidable.

\subsection*{Statements of Results}

The main purpose of the present note is to prove the following result:

\begin{thm}
\label{mainthm}
For each $d\geq 2$ there exists an oracle Turing machine $M^\phi$, with the oracle
representing the coefficients of a polynomial of degree $d$ with the following 
property.
Suppose that $p(z)$ is a polynomial of degree $d$ such that the filled Julia set $K_p$ has 
no interior (so that $J_p=K_p$). Then $K_p$ is computable by $M^\phi$.
\end{thm}

\noindent
The known results in this direction, which make the statement of the theorem interesting, are the 
following. Independently, the second author \cite{thesis} and Rettinger \cite{Ret} have demonstrated that
{\it hyperbolic} Julia sets are polynomial-time computable by an oracle TM.
In sharp contrast, in \cite{BY} the second and third authors demonstrated the existence of non-computable
Julia sets in the family of quadratic polynomials. The latter examples are given by polynomials
with Siegel disks. It is shown in \cite{BY} that there exist non-computable examples with rather
wild topology. By a method of Buff and Ch{\'e}ritat \cite{BC} the boundary of the Siegel disk 
can be made smooth, and as a consequence, the critical point is not in the boundary and 
the Julia set is not locally connected (see \cite{Ro} for a discussion of the
topological anomalies in such Julia sets).

As our present theorem shows, however, the notions of topological complexity and computational complexity
are rather distinct. It covers, for example, the case of Cremer quadratics, which are also known
to possess non-locally connected Julia sets. Further, we find the following statement plausible:

\medskip
\noindent
{\bf Conjecture.}
{\it The procedure of \cite{BY} can be carried out so as to yield non-computable and locally connected
Siegel quadratics.}

\medskip

\noindent
Our main result implies that the reason one is not able to produce high resolution computer-generated
images of Cremer Julia sets is due to the enormous computational time required
(whether such images would be very informative is another question).
 In view of this, it is of great interest to obtain a classification of Julia
sets into classes of computational hardness.

\medskip
\noindent
Finally, let us state a straightforward generalization of the main theorem
to the case of a rational map without rotation domains (either Siegel disks or Herman rings).
First note that by the Fatou-Sullivan Classification, every Fatou component of such a rational
map belongs to the basin of a (super)attracting or parabolic orbit. 
A rational map $f$ of degree
$d$ can have at most $2d-2$ of such orbits by Shishikura's result. 
Of these, the attracting orbits can be found algorithmically, with an arbitrary precision,
as the condition $|(f^m)'(x)|<1$ can be verified with a certainty by a TM.
For an algorithm to ascertain the
position of a parabolic orbit a finite amount of information is 
needed. For example, we can ask for its period $m$,
and a good rational approximation to distinguish the point among the roots of $f^m(z)=z$.

\begin{thm}
\label{secondthm}
Let $f$ be a rational map $f:\hat \CC\to\hat\CC$ without rotation domains.
Then its Julia set is computable in the spherical metric by an oracle Turing machine $M^\phi$ with 
the oracle representing the coefficients of $f$. The algorithm uses the nonuniform
information on the number and positions of the parabolic orbits,
as described above, as well as the (rational) values of the $\frac{1}{2 \pi} \log$'s of the multipliers of the parabolic
orbits.
\end{thm}

\noindent
In particular, if there are no parabolic orbits the computation can be performed 
by a single oracle machine. 
As a function computed by an oracle Turing machine must be continuous on its domain (see e.g. \cite{thesis}),
we have the following corollary:

\begin{cor}
Denote $\cR^d\subset \CC^{2 d-2}$ the set of parameters of degree $d$ rational maps  without
 rotation domains or parabolic orbits. 
Then a Julia set depends continuously, in Hausdorff sense on $\hat \CC$, on the parameter in $\cR^d$.
\end{cor}

\noindent
While this is well-known (see the discussion in \cite{Do}), it is notable that the proof
of this dynamical statement is produced by a computer science argument. \footnote{We owe this
comment to J. Milnor}

 \medskip

\noindent
{\bf Acknowledgement.} The authors wish to thank John Milnor, whose encouragement and
questions have inspired this work.

\section{Proof of the Main Theorem}

\subsection*{Idea of the proof} The idea of the argument is to find two sequences of Hausdorff bounds for the
filled Julia set. The {\it inside} ones are given by rational approximation of the periodic cycles
of $p(z)$, accomplished by approximately solving a sequence of polynomial equations. 
The {\it outside} ones are given by preimages of a large disk centered at the origin. The algorithm stops
when the difference between the two bounds is sufficiently small. This will always happen, provided
the interior of the filled Julia set is empty.

We make no claims about the efficiency of our algorithm. Indeed, in most cases it will be far
too slow to be practical. The existence of an  efficient algorithm for the challenging
cases, such as Julia sets with Cremer points, is an interesting open problem.

\begin{lem}
\label{lb}
For every natural $n$ we can compute a sequence of rationals $\{q_i\}$ 
such that 
\beq
\label{lbeq}
B(J_p,2^{-(n+2)})  \Subset \bigcup_{i=1}^{\infty} B(q_i, 2^{-(n+1)}) 
\Subset B(J_p, 2^{-n}).
\eeq
\end{lem}

\begin{proof}
First of all, note that for any polynomial $Q(z)$, we can list rational approximations 
$r_1, r_2, \ldots, r_m$ of all the roots 
$\al_1, \al_2, \ldots, \al_m$ of $Q(z)-z$ with an arbitrarily good precision $2^{-k}$ (a classical
reference is \cite{Wey}). 

Let $M>0$ be some bound on $|Q''(z)|$ in the area of the roots. Then $|Q'(r_i)|>1+2^{-k}M$ 
implies that $|Q'(\al_i)|>1$, and in fact acts as a {\em certificate} that this is the case. 
It is easy to see that if, in fact, $|Q'(\al_i)|>1$, then eventually we will find a rational
approximation $r_i$ of $\al_i$ for which we {\em know} that $|Q'(\al_i)|>1$. 

We now use the fact that the Julia set $J_p$ is the closure of the 
repelling periodic orbits of $p$. 
Thus we can compute a sequence $\{q_i\}$ of rational approximations of all the repelling periodic 
points such that:
\begin{itemize}
\item 
For each $i$, there is a repelling periodic point $\al_i$ with 
$|q_i-\al_i|<2^{-(n+3)}$, and
\item 
for each repelling periodic point $\al$, a $2^{-(n+3)}$-approximation of 
$\al$ will 
eventually appear in $\{q_i\}$. 
\end{itemize}

We claim that \eref{lbeq} holds for the sequence $\{q_i\}$.

For any $x \in B(J_p,2^{-(n+2)})$ there is a $y\in J_p$ with 
$|x-y|<2^{-(n+2)}$. There is a repelling periodic point $z$ with 
$|z-y|<2^{-(n+3)}$, and the rational point $q_i$ approximating it 
satisfies 
$|z-q_i|<2^{-(n+3)}$, hence 
$$
|x-q_i| \leq |x-y|+|y-z|+|z-q_i| < 2^{-(n+2)} + 2^{-(n+3)} +2^{-(n+3)} =
2^{-(n+1)}.
$$
So $x \in B(q_i, 2^{-(n+1)})$. 

On the other hand, for any $q_i$ there is a repelling $z\in J_p$ with 
$|q_i-z|<2^{-(n+3)}$, hence $B(q_i, 2^{-(n+1)}) \subset B(J_p, 
2^{-(n+1)}+2^{-(n+3)}) \subset B(J_p, 2^{-n})$. 
\end{proof}

\begin{proof}[Proof of \thmref{mainthm}] 
Fix $m\in\NN$. Our algorithm to find $C_m\in\cC$ works as follows.
Evaluate a large number $b>0$ such that if $|z|>b$ then $|p(z)|>b$ and the orbit $p^n(z)\to\infty$.
Set $D\equiv B(0,b)$.
At the $k$-th step:

\begin{itemize}
\item  compute the finite union $B_k=\cup_{i=1}^k B(q_i,2^{-(m+1)})\in \cC$ as in the previous
lemma; 
\item compute with precision $2^{-(m+3)}$ the preimage $p^{-k}(D)$, that is,
find $D_k\in \cC$ such that 
$$d_H(\overline{D_k},\overline{p^{-k}(D)})<2^{-(m+3)};$$
\item if $D_k\subset B_k$ output $C_m=B_k$ and terminate. Otherwise, go to step $k+1$.
\end{itemize}

\noindent
We claim that this algorithm eventually terminates, and outputs an approximation $C_m\in\cC$ with
$d_H(C_m,K_p)<2^{-m}$. 

\noindent
{\bf Termination.} First, we note that  $p^{-k}(D)$ is compactly contained in $p^{-(k-1)}(D)$
for all $k\in\NN$. Since in our case $J_p=K_p$, we have
\begin{equation}
\label{K equal}
K_p=\bigcap_{k\in \NN}p^{-k}(D)=\overline{\{\text{repelling periodic orbits of }p\}}.
\end{equation}
By the previous lemma, there exists $k_0\in\NN$ such that for all $k\geq k_0$ we have
$$B(K_p,2^{-(m+2)})\Subset B_k\Subset B(K_p,2^{-m}).$$
By (\ref{K equal}) there exists $k_1\in\NN$ such that for all $k\geq k_1$ we have
$$p^{-k}(D)\subset B(K_p,2^{-(m+3)}),\text{ and hence }D_k\subset B(K_p,2^{-(m+2)}).$$
Therefore, for $k\geq\max(k_0,k_1)$ we have $D_k\subset B_k$ and the algorithm will terminate.

\medskip
\noindent
{\bf Correctness.} Now suppose that the algorithm terminates on step $k$. Since $D_k\subset B_k$
and $K_p\subset B(D_k,2^{-(m+3)})$ we have $K_p\subset B(C_m,2^{-(m+3)})$. On the other hand,
$\cup \{ q_i\}\subset K_p$, and thus $B_k=C_m\subset B(K_p,2^{-(m+1)})$.

\end{proof}

\noindent
\begin{proof}[Proof of \thmref{secondthm}]
We will have to modify the way the set $D_k$ is generated in the preceding algorithm.
Note that the procedure which finds periodic points of $f$ of period $m$, and 
verifies the conditions $|(f^m)'(x)|>1$, or $|(f^m)'(x)|<1$ with precision $2^{-m}$
will {\it eventually} certify every attracting periodic point.
For each such point, we will also identify a Euclidean disk
$U_x\in\cC$ containing it
with the property $f^m(U_x)\Subset U_x$. 

 Further, for each parabolic point $x$ of period $l$ let 
$U^k_x\in\cC$ be a topological disk compactly contained in the immediate 
basin of attraction of $x$ (i.e. the union of attracting petals), with the properties:
$$U^k_x\subset U^{k+1}_x\text{, and }\cup_k U^k_x\text{ is an attracting petal}.$$
Such neighborhoods can easily be produced algorithmically given the non-uniform information
on the parabolic orbits as stated in the Theorem.
Our algorithm will then work as before, with the difference that the set
$D_k$ is defined as a $2^{-(m+3)}$-approximation of the {\it complement} of the preimage
$$W_k=f^{-k}(\bigcup_{x\text{ is attracting located so far}} U_x )\cup f^{-k}(\bigcup_{x\text{ is parabolic}}
 U^k_x).$$
By the Sullivan's Non-Wandering Theorem, we have $J_p=\cap(\hat\CC\setminus W_k)$. The 
algorithm is verified exactly as in the previous theorem.
\end{proof}

\section{Conclusion.}

We have shown that only rational maps  with rotation domains (Siegel disks or Herman rings)
may have uncomputable Julia sets. As was demonstrated earlier by the second
and third authors \cite{BY}, such a possibility is indeed realized. 
To understand the cause of this phenomenon, consider the
practical problem of plotting the Julia set of a Siegel quadratic $P_\theta(z)=z^2+e^{2\pi i\theta}z$.
The preimages of a large disk around the Julia set will rapidly converge to $K_{P_\theta}$. However,
when they are already quite close to $J_{P_\theta}$, these preimages will still be quite far
from the Siegel disk $\Delta_\theta$, as they will not fill the narrow fjords leading towards
$\partial \Delta_\theta$. 
To have a faithful picture of a $2^{-n}$-neighborhood of $J_{P_\theta}$, one then needs to 
draw a contour inside the Siegel disk, close to its boundary. In effect, this means knowing
the inner radius of the Siegel disk.
Indeed, we have the following easy statement:

\begin{prop}
\label{comp siegel}
Suppose the inner radius $\rho_\theta$ of $\Delta_\theta$ is known. 
Then the Julia set $J_{P_\theta}$ is computable.
\end{prop}

\begin{pf}
The algorithm to produce the $2^{-n}$ approximation of the Julia set is the following.
First, compute a large disk $D$ around $0$ with $P_\theta(D)\Supset D$. Then,
\begin{itemize}
\item[(I)] compute a set $D_k\in\cC$ which is a 
$2^{-(n+3)}$-approximation of the preimage $P_\theta^{-k}(D)$;
\item[(II)] set $W_k$ to be the round disk with radius $\rho_\theta-2^{-k}$ about the
origin. Compute a set $B_k\in \cC$ which 
is a  $2^{-(n+3)}$-approximation of $P_\theta^{-k}(W_k)$;
\item[(III)] if $D_k$ is contained in a $2^{-(n+1)}$-neighborhood of $B_k$, then output
a $2^{-(n+1)}$-neighborhood of $D_k - B_k$, and stop. If not, go to step (I).
\end{itemize} 

\noindent
A proof of the validity of the algorithm follows along the same lines as that in
\thmref{mainthm}.

\end{pf}

\noindent
The well-recognizable images of the ``good'' Siegel Julia sets, such as the one
for $\theta_*=(\sqrt{5}-1)/2$, the golden mean, are indeed computable by this recipe:

\begin{prop}
The inner radius $\rho_{\theta_*}$ is a computable real number.
\end{prop}
\begin{pf}
We appeal to the renormalization theory for golden-mean Siegel disks (see \cite{McM}),
which implies, in particular, that the boundary of $\Delta_{\theta_*}$ is 
self-similar up to an exponentially small error. As a consequence, 
there exist $C>0$, and $\lambda>1$ such that 
$$\inf\{|f^i(0)|,\; i=0,\ldots,q_n\}<\rho_{\theta_*}(1+C\lambda^{-n}),$$
where $p_n/q_n$ is the $n$-th continued fraction convergent of $\theta_*$.
\end{pf}

\begin{figure}[h]
\centerline{\psfig{figure=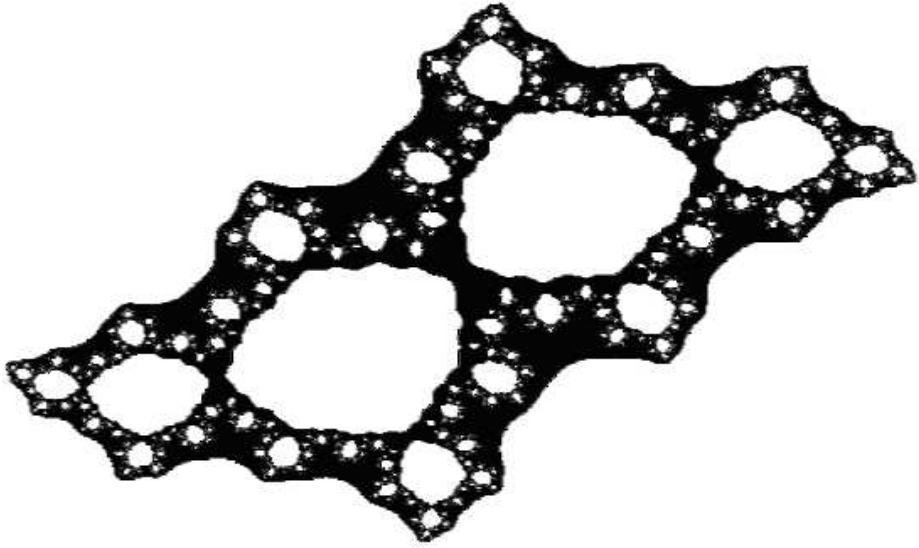,width=12.5cm}}
\caption{An image of the Siegel Julia set $J_{\theta_*}$ produced by the algorithm of 
\propref{comp siegel}}
\end{figure}

\noindent
In contrast, in the non-computable examples constructed in \cite{BY} the 
boundary of the Siegel disk oscillates {\it unpredictably} at infinitely many scales.
More precisely, recall the definition of the {\it conformal radius} $r(U,u)$ of a
simply-connected domain $U$ with a marked point $u$ in its interior. For the unique
holomorphic bijection $\Psi$ from the unit disk $\DD$ to $U$ with $\Psi(0)=u$, and
$\Psi'(0)>0$ set 
$$r(U,u)=\Psi'(0).$$
By the Koebe 1/4 Theorem of classical complex analysis, the inner radius of $U$ as seen
from $u$ is at least $\frac{1}{4}r(U,u)$.

Now set $r_\theta=r(\Delta_\theta,0)$. It is shown in \cite{BY} that there exist parameter
values $\tau\in\RR/\ZZ$ for which the function $\theta\mapsto r_\theta$ is uncomputable
on the set consisting of a single point $\{\tau\}$. On the other hand, we have:

\begin{prop}
\label{pr1}
Computability of $J_{P_\theta}$ implies computability of $r_\theta$ for Siegel parameter values $\theta$.
\end{prop}

\noindent
We note that a converse is also true, so computing the conformal radius is  the only
possible obstruction to drawing $J_{P_\theta}$:

\begin{prop}
\label{pr2}
Suppose $r_\theta$ is computable, then so is $J_{P_\theta}$.
\end{prop}

\noindent
We will require the following facts of complex analysis:

\begin{lem}
\label{lem1}
Let $(U,u)$ be as above, and suppose $U_1$ is a simply-connected subdomain of $U$
containing the point $u$. Let $\eps>0$ be such that $\partial U_1\subset B(\partial U,\eps)$.
Denote $\rho$ and $\rho_1$ the inner radii of $U$ and $U_1$ respectively, as seen from $u$.
 There exist explicit positive constants $C_1=C_1(r(U,u))$, and  $C_2=C_2(r(U,u))$ such that
$$ C_1|\rho-\rho_1|\leq |r(U,u)-r(U_1,u)|\leq C_2\sqrt{\eps}.$$

\end{lem}

\noindent
The second inequality is an elementary consequence of Koebe Theorem and the
 constant may be chosen as $C_2=4 \sqrt{r(U,u)}$ (see \cite{RZ}, for example). 
The other inequality may be found in many advanced texts
on Complex Analysis, e.g. \cite{Pom}.

\begin{proof}[Proof of \propref{pr2}]
We will show the computability of the inner radius of $\Delta_\theta$. The algorithm works as 
follows:
\begin{itemize}
\item[(I)] For $k\in\NN$  compute a set $D_k\in\cC$ which is a 
$2^{-m}$-approximation of the preimage $P_\theta^{-k}(D)$, for some sufficiently
large disk $D$;
\item[(II)]  evaluate the
conformal radius $r(D_k,0)$ with precision $2^{-(m+1)}$ (this can be done, for example,
by using one of the numerous existing methods for computing the Riemann Mapping of a computable domain,
see \cite{Zhou}, for example), if necessairy, compute $D_k$ or parts of $D_k$ with 
a higher degree of precision;
\item[(III)] if $r(D_k,0)$ is $2^{-m}{C_1}$-close to $r_\theta$ then compute the inner radius
$\rho(D_k)$ around $0$ with precision $2^{-m}$ and output this number. Else, increment $k$ and
return to step (I).
\end{itemize}

\medskip
\noindent
{\bf Termination.}
Let $K=K_{P_\theta}$ be the filled Julia set of $P_\theta$. Then $$\cap_{k=0}^{\infty} D_k = K \supset \Delta_\theta$$
and $D_0 \supset D_1 \supset D_2 \supset \ldots$. Hence there will be a step $k_0$ after which
$$\dist(\partial D_k,J_{P_\theta})<\frac{2^{-2 m}C_1^2}{C_2^2}.$$
$\partial \Delta_\theta \subset J_{P_\theta}$, and by \lemref{lem1} this implies that 
$$|r(D_k,0)-r(\Delta_\theta,0)|=|r(D_k,0)-r_\theta|<2^{-m} \frac{C_1}{C_2} C_2 = 2^{-m} C_1,$$
and the algorithm will terminate on step (III).

\medskip
\noindent
{\bf Correctness.}
Now suppose the algorithm has terminated on step (III). $\Delta_\theta \subset 
D_k$, and \lemref{lem1} implies that 
$$|\rho(D_k)-\rho_\theta|\leq 2^{-m}.$$

\end{proof}

\comm{\begin{proof}[Proof of \propref{pr2}]
We will show the computability of the inner radius of $\Delta_\theta$. The algorithm works as 
follows:
\begin{itemize}
\item[(I)] For $k\in\NN$ compute the finite union $B_k=\cup_{i=1}^k B(q_i,2^{-(2 m+2)})\in \cC$
as in \lemref{lbeq};
\item[(II)] for the connected component $\Omega_k$ of $\CC\setminus B_k$ containing $0$ evaluate the
conformal radius $r(\Omega_k,0)$ with precision $2^{-(m+1)}$ (this can be done, for example,
by using one of the numerous existing methods for computing the Riemann Mapping of a domain with
an explicitly described boundary);
\item[(III)] if $r(\Omega_k,0)$ is $2^{-(m-1)}C_2$-close to $r_\theta$ then compute the inner radius
$\rho(\Omega_k)$ around $0$ with precision $2^{-m}$ and output this number. Else, increment $k$ and
return to step (I).
\end{itemize}

\medskip
\noindent
{\bf Termination.}
Given the density of the points $\{q_i\}$ in $J_p$, there will be a step $k_0$ after which
$$\dist(\partial\Omega_k,\Delta_\theta)<2^{-2 m}.$$
By \lemref{lem1} this implies that 
$$|r(\Omega_k,0)-r_\theta|<2^{-m}C_2,$$
and the algorithm will terminate on step (III).

\medskip
\noindent
{\bf Correctness.}
Now suppose the algorithm has terminated on step (III). \lemref{lem1} implies that 
$$|\rho(\Omega_k)-\rho_\theta|\leq 2^{-(m-1)}\frac{C_2}{C_1}.$$

\end{proof}
} 

\end{document}